\documentclass[preprint,12pt]{elsart}
\usepackage{amssymb,amsmath,graphicx}
\usepackage{longtable}
\begin{document}
\begin{frontmatter}
\title{A degenerate kernel method for eigenvalue problems of a class of non-compact operators}
\author[label1]{H.~Majidian}
\ead{majid@modares.ac.ir}
\author[label2]{E.~Babolian\corauthref{c1}}
\ead{babolian@saba.tmu.ac.ir}
\address[label1]{Department of Mathematics, Tarbiat Modares University, PO~Box 14115-175 Tehran, Iran}
\address[label2]{Department of Mathematics, Tarbiat Moallem University, Tehran 1561836314, Iran}
\corauth[c1]{Corresponding author. Tel: +98 21 7750 7722   Fax:
+98 21 7760 2988}
\begin{abstract}
We consider the eigenvalue problem of certain kind of non-compact
linear operators given as the sum of a multiplication and a kernel
operator. A degenerate kernel method is used to approximate
isolated eigenvalues. It is shown that entries of the
corresponding matrix of this method can be evaluated exactly. The
convergence of the method is proved; it is proved that the
convergence rate is $O(h)$. By some numerical examples, we
confirm the results.
\end{abstract}
\end{frontmatter}
\section{Introduction}
Consider the following non-compact operator defined on
$X=\mathcal{L}^{\infty}[a,b]$ by
\begin{equation}\label{operator}
\mathcal{A}\,f(x)=\alpha\int_{a}^{b}k(x,u)~f(u)~du-x^{2}~f(x),\qquad
\forall x\in[a,b].
\end{equation}
Here, $[a,b]$ is a compact interval,
$k(.,.)\in\mathcal{C}^{0}([a,b]\times [a,b])$, and $\alpha$ is a
real constant. The non-compactness of the operator $\mathcal{A}$
is due to the involved multiplication operator with the
corresponding function $x^{2}$ \cite{morrison}.

We are concerned with the numerical approximation of isolated
eigenvalues of the following eigenvalue problem:
\begin{equation}\label{eigprob}
\lambda\,f=\mathcal{A}\,f.
\end{equation}

Eigenproblem~(\ref{eigprob}) has some applications in
electromagnetism \cite{shahabadi}. Also, a more general form
problem~(\ref{eigprob}) describes the COA model in population
genetics \cite{burger}. In \cite{redner}, this model has been
treated numerically by Nystr\"{o}m and Galerkin methods in the
presence of some extra assumptions. A degenerate kernel method
with piecewise linear interpolation with respect to the variable
$u$ has been developed for problem~(\ref{eigprob}) in
\cite{majidian1}. In this method, isolated eigenvalues of
problem~(\ref{eigprob}) are approximated by eigenvalues of an
operator $\mathcal{A}_{n}$ with rank $n$; but the entries of the
matrix eigenvalue problem associated with $\mathcal{A}_{n}$ can
not be evaluated exactly. Hence, we need to replace various
integrals by numerical quadratures. This perturbation in entries
causes unacceptable approximations for ill-conditioned
eigenvalues.

In \cite{kaneko} a degenerate kernel method is proposed for the
numerical solutions of integral equations of the second kind with
a smooth kernel which does not need any numerical quadrature, and
the entries of its corresponding matrix are evaluated exactly if
we ignore round-off errors. Genaneshwar \cite{gnan} has been
extended this method to the corresponding eigenvalue problem with
the same benefit. In this paper we extend this method for the
numerical solutions of eigenproblem~(\ref{eigprob}).

The rest of this paper is organized as follows: In
Section~\ref{secmethod}, we describe the degenerate kernel method
for numerical solutions of problem~(\ref{eigprob}). Convergence of
the method and its rate are discussed in Section~\ref{secmethod}.
In Section~\ref{secexmp}, we give some numerical results to
illustrate the accuracy of the method.

We recall that $X=\mathcal{L}^{\infty}[a,b]$ throughout this
paper.
\section{Degenerate kernel method}\label{secmethod}
We begin this section by some preliminaries.

Let $Y$ be any Banach space over the complex field $\mathbb{C}$.
We denote the space of bounded linear operators from $Y$ into $Y$
by $BL(Y)$. Let $T\in BL(Y)$. The \emph{resolvent set} of $T$ is
given by
$$\rho(T)=\{z\in\mathbb{C}\,:\,(T-zI)^{-1}\in BL(Y)\}.$$
The \emph{spectrum} of $T$, denoted by $\sigma(T)$, is defined as
$\sigma(T)=\mathbb{C}\backslash\rho(T)$. The \emph{point
spectrum} of $T$ consists of all $\lambda\in\sigma(T)$ such that
$T-\lambda I$ is not one-to-one. In this case $\lambda$ is called
an \emph{eigenvalue} of $T$. If $\lambda$ is an eigenvalue of
$T$, then the smallest positive integer $l$ such that
$ker(T-\lambda I)^{l}=ker(T-\lambda I)^{l+1}$ is called the
\emph{ascent} of $\lambda$. The dimensions of $ker(T-\lambda I)$
and $ker(T-\lambda I)^{l}$ are called \emph{geometric
multiplicity} and \emph{algebraic multiplicity} of $\lambda$,
respectively. If the algebraic multiplicity of $\lambda$ equals
one, then it is called a \emph{simple eigenvalue} of $T$.

We now describe the degenerate kernel method. Define the
operators $\mathcal{K}$ and $\mathcal{M}$ on $X$ as follows:
\begin{equation}
\mathcal{K}\,f(x)=\alpha\int_{a}^{b}\,k(x,u)\,f(u)\,du,
\end{equation}
and
\begin{equation}
\mathcal{M}\,f(x)=x^{2}\,f(x).
\end{equation}
Then $\mathcal{A}=\mathcal{K}-\mathcal{M}$.

For an integer $n>1$, consider the following partition of $[a,b]$:
\begin{equation*}
a=x_{0}<x_{1}<\dots<x_{n-1}<x_{n}=b.
\end{equation*}
Let $I_{j}=[x_{j-1},x_{j}]$, $h_{j}=x_{j}-x_{j-1}$ for
$j=1,2,\dots,n$. Also, let $h=\max\limits_{1\le j\le n}h_{j}$
denotes the norm of the partition. We assume that $h\to0$ as
$n\to\infty$. For a positive integer $r$, let
$B_{r}=\{\tau_{1},\dots,\tau_{r}\}$ be the set of $r$ Gauss
points, i.e., the zeros of the Legendre polynomial
$\frac{d^{r}}{dt^{r}}(t^{2}-1)^{r}$ in the interval $[-1,1]$.
Define $f_{j}:[-1,1]\to I{j}$ as follows:
\begin{equation*}
f_{j}(t)=\frac{\displaystyle 1-t}{\displaystyle 2}x_{j-1}+
\frac{\displaystyle 1+t}{\displaystyle 2}x_{j},\quad t\in[-1,1].
\end{equation*}
Then
$S=\cup_{j=1}^{n}\,f_{j}(B_{r})=\{\xi_{ij}=f_{j}(\tau_{i})\,:\,i=1,\dots,r,\,j=1,2,\dots,
n\}$ is the set of $N_{h}=nr$ Gauss points on $[a,b]$. For each
$i=1,\dots,r$, let
\begin{equation*}
\ell_{i}(x)=\frac{\displaystyle
(x-\tau_{1})\dots(x-\tau_{i-1})(x-\tau_{i+1})\dots(x-\tau_{r})}{\displaystyle
(\tau_{i}-\tau_{1})\dots(\tau_{i}-\tau_{i-1})(\tau_{i}-\tau_{i+1})\dots(\tau_{i}-\tau_{r})},\quad
x\in [-1,1],
\end{equation*}
be the Lagrange polynomial of degree $r-1$ on $[-1,1]$, which
satisfies $\ell_{i}(\tau_{j})=\delta_{ij}$. It is notable that
for $r=1$, $\ell_{1}(x)=1$ on $[-1,1]$. Define
\begin{equation*}
  \phi_{pq}(x)=
  \begin{cases}
    \ell_{p}(f^{-1}_{q}(x)), & x\in I_{q}, \\
    0, & \text{otherwise}.
  \end{cases}
\end{equation*}
Then $\phi_{pq}(\xi_{ij})=\delta_{pi}\delta_{jq}$, for
$p,i=1,\dots,r$ and $j,q=1,2,\dots,n$.

We set the following notations $t_{(j-1)r+i}=\xi_{ij}$,
$\psi_{(j-1)r+i}=\phi_{ij}$, for $i=1,\dots,r$ and
$j=1,2,\dots,n$. Then $S=\{t_{i}\,:\,i=1,2,\dots,N_{h}\}$ is the
set of $nr$ Gauss points in $[a,b]$. Now define the degenerate
kernel by
\begin{equation*}
k_{N_{h}}(x,u)=\sum\limits_{i=1}^{N_{h}}\sum\limits_{j=1}^{N_{h}}\,k(t_{i},t_{j})
\psi_{i}(x)\psi_{j}(u),\qquad x,u\in[a,b].
\end{equation*}
The approximated operator $\mathcal{K}_{N_{h}}$ is defined by
\begin{equation}
\mathcal{K}_{N_{h}}\,f(x)=\alpha\int_{a}^{b}\,k_{N_{h}}(x,u)\,f(u)\,du.
\end{equation}
For the approximation of the multiplication operator
$\mathcal{M}$, define the piecewise constant functions
$s_{N_{h}}(x)$ on $[a,b]$ as follows:
\begin{equation*}
\text{If}\;x\in I_{i},\;\text{then}\;s_{N_{h}}(x)=m_{i}=
\left(\frac{\displaystyle x_{i-1}+x_{i}}{\displaystyle
2}\right)^{2}, \qquad i=1,2,\dots,N_{h}.
\end{equation*}

Now define the approximate operators $\mathcal{K}_{N_{h}}$ and
$\mathcal{M}_{N_{h}}$ on $X$ as follows:
\begin{equation}
\mathcal{M}_{N_{h}}\,f(x)=s_{N_{h}}(x)\,f(x).
\end{equation}
Consider now the following approximate eigenvalue problem on $X$:
\begin{equation}\label{approxeigprob}
\lambda_{N_{h}}\,f_{N_{h}}(x)=\mathcal{A}_{N_{h}}\,f_{N_{h}}(x),
\end{equation}
where
$\mathcal{A}_{N_{h}}=\mathcal{K}_{N_{h}}-\mathcal{M}_{N_{h}}$.
The rank of the operator $\mathcal{A}_{N_{h}}$ is finite, so the
eigenvalue problem~(\ref{approxeigprob}) is equivalent to a matrix
eigenvalue problem as we show in the following:

Eigenproblem~(\ref{approxeigprob}) expands into
\begin{equation}\label{approxeigprobexpand}
\lambda_{N_{h}}\,f_{N_{h}}(x)=
\alpha\,\sum\limits_{p=1}^{N_{h}}\sum\limits_{j=1}^{N_{h}}\,
k(t_{p},t_{j})\,\psi_{p}(x)\,c_{j}-s_{N_{h}}(x)\,f_{N_{h}}(x),
\end{equation}
where
\begin{equation*}
c_{j}=\int_{a}^{b}\psi_{j}(u)f_{N_{h}}(u)\,du.
\end{equation*}
For each $i=1,2,\dots,N_{h}$, we multiply both sides of
(\ref{approxeigprobexpand}) by $\psi_{i}(x)$ and integrate over
the interval $[a,b]$ with respect to the variable $x$ to obtain
\begin{equation}
\lambda_{N_{h}}\,c_{i}=\alpha\,\sum\limits_{j=1}^{N_{h}}k_{i,j}\,c_{j}-m_{i}\,c_{i},
\qquad i=0,1,\dots,n,
\end{equation}
where
\begin{equation*}
k_{i,j}=\sum\limits_{p=1}^{N_{h}}k(t_{p},t_{j})\int_{a}^{b}\psi_{i}(x)\,\psi_{p}(x)\,dx.
\end{equation*}
Therefore, eigenvalues of the approximate
problem~(\ref{approxeigprob}) are the eigenvalues of the matrix
$A=(a_{i,j})$ of order $N_{h}$, defined as
\begin{equation}\label{matrixa}
a_{i,j}=\alpha\,k_{i,j}-m_{i}\,\delta_{i,j},\qquad
i,j=1,2,\dots,N_{h}.
\end{equation}
It is easy to see that
%The entries of $A$ depends on the values
$$
\int_{a}^{b}\psi_{i}(x)\,\psi_{p}(x)\,dx=0,\quad \text{if}\, i\ne
p.
$$
Therefore, $k_{ij}$ reduces to
\begin{equation}
k_{i,j}=k(t_{i},t_{j})\,\int_{a}^{b}\psi_{i}^{2}(x)\,dx.
\end{equation}

Since $\psi_{i}$ for $i=1,2,\dots,N_{h}$ are basis functions
having small support on $[a,b]$ with the simple structure, the
integrals $\int_{a}^{b}\psi_{i}^{2}(x)\,dx$ can be evaluated
exactly. Values of these integrals depend on the length of their
support. Thus, if $x_{i}$ for $i=1,2,\dots,n$ are equidistance
points and the kernel $k(.,.)$ is Hermitian, then the matrix $A$
is normal and as a result, each eigenvalue of $A$ is
well-conditioned (see \cite{saad}).
\section{Convergence and the rate of convergence}\label{secconv}
\begin{thm}\label{thrmconverg}\cite{ahues}
Let $Y$ be a Banach space. Let also $T$ and
$\{T_{n}\}_{0}^{\infty}$ be operators in $BL(Y)$. Assume that the
sequence $\{T_{n}\}$ converges in norm to $T$, i.e.,
$\|T_{n}~-~T\|\to0$ as $n\to\infty$. Let $\lambda$ be an isolated
point of $\sigma(T)$. For each positive
$\epsilon<dist(\lambda,\sigma(T)\backslash\{\lambda\})$, define
$$\Lambda_{n}:=\{\lambda_{n}\in\sigma(T_{n})\,:\,|\lambda_{n}-\lambda|<\epsilon\}.$$
Then for $n$ large enough, $\Lambda_{n}\neq\emptyset$ and if
$\lambda_{n}\in\Lambda_{n}$, the sequence $\{\lambda_{n}\}$
converges to $\lambda$.
\end{thm}
In order to show that the degenerate kernel method proposed in
previous section is convergent, it is enough to show that the
operators $\mathcal{A}$ and $\{\mathcal{A}_{N_{h}}\}$ satisfy
conditions of Theorem~\ref{thrmconverg}, i.e., they are in $BL(X)$
and $\mathcal{A}_{N_{h}}$ converges in norm to $\mathcal{A}$ as
$N_{h}\to\infty$. The linearity of these operators is obvious. In
\cite{majidian1}, we have shown that the operator $\mathcal{A}$ is
bounded. By a similar discussion it is seen that the operator
$\mathcal{A}_{N_{h}}$ is bounded, too.

If we assume that $k(.,.)\in\mathcal{C}^{r}([a,b]\times [a,b])$,
then $\|\mathcal{K}-\mathcal{K}_{N_{h}}\|=O(h^{r})$ \cite{gnan}.
Also, from \cite{majidian1} we have
$\|\mathcal{M}-\mathcal{M}_{N_{h}}\|=O(h)$. Therefore,
\begin{equation}\label{oa}
\|\mathcal{A}-\mathcal{A}_{N_{h}}\|=O(h).
\end{equation}
On the other hand, with the assumptions of
Theorem~\ref{thrmconverg}, if $\lambda$ is a simple eigenvalue of
$\mathcal{A}$, then there exists a constant $c$ such that
\cite[Page 201]{ahues}
\begin{equation}\label{neqeig}
|\lambda_{N_{h}}-\lambda|\le
c\,\|\mathcal{A}_{N_{h}}-\mathcal{A}\|,\qquad\text{for all large
$N_{h}$}.
\end{equation}
From (\ref{oa}) and (\ref{neqeig}) we conclude that the proposed
degenerate kernel method of previous section is convergent with
the rate $O(h)$ for simple eigenvalues.
\section{Numerical results}\label{secexmp}
In this section, we solve two samples of problem~(\ref{eigprob})
using the degenerate kernel method, proposed in
Section~\ref{secmethod}. In order to show the convergence, the
problems are solved for different numbers of interpolation nodes.
Two arbitrary positive eigenvalues of this problem, denoted by
$\lambda(1)$ and $\lambda(2)$, are approximated by
$\lambda_{n}(1)$ and $\lambda_{n}(2)$, respectively. In both
examples, we take $r=2$ and compute the absolute error and the
ratio of each eigenvalue.
\subsection{Example~1}
Consider problem~(\ref{eigprob}) with the kernel
$$k(x,u)=\exp(-(u-x)^{2}),$$ and the constant $\alpha=1$ on the interval
$[-2,2]$. Computational results are shown in Table~\ref{tblerr1}.
\begin{table}
\caption{Computational results for sample
problem~1.}\label{tblerr1}
\begin{tabular}{lllll}
\hline $n$ & $|\lambda_{n}(1)-\lambda(1)|$ & $ratio(1)$ &
$|\lambda_{n}(2)-\lambda(2)|$ & $ratio(2)$ \\
\hline
10  & 0.0763 &       & 0.0916 &      \\
20  & 0.0398 & 1.92  & 0.0460 & 1.99 \\
40  & 0.0202 & 1.97  & 0.0229 & 2.01 \\
80  & 0.0100 & 2.02  & 0.0113 & 2.03 \\
160 & 0.0048 & 2.08  & 0.0054 & 2.09 \\ \hline
\end{tabular}
\end{table}
\subsection{Example~2}
Consider problem~(\ref{eigprob}) with the kernel
$$
    k(x,u)=
    \frac{\displaystyle 1}{\displaystyle 1+(u-x)^{2}}
$$
and the constant $\alpha=1$ on the interval $[-4,4]$.
Computational results are shown in Table~\ref{tblerr2}.

In both examples, convergence of solutions is seen. The ratio of
error in each step also confirms our theoretical discussion in
Section~\ref{secconv}.
\begin{table}
\caption{Computational results for sample
problem~2.}\label{tblerr2}
\begin{tabular}{lllll}
\hline $n$ & $|\lambda_{n}(1)-\lambda(1)|$ & $ratio(1)$ &
$|\lambda_{n}(2)-\lambda(2)|$ & $ratio(2)$ \\
\hline
10  & 0.1607 &       & 0.0712 &      \\
20  & 0.0793 & 2.03  & 0.0287 & 2.48 \\
40  & 0.0396 & 2.00  & 0.0124 & 2.31 \\
80  & 0.0196 & 2.02  & 0.0055 & 2.25 \\
160 & 0.0095 & 2.06  & 0.0025 & 2.20 \\ \hline
\end{tabular}
\end{table}
\section{Acknowledgements} The authors wish to thank Masoud Amini for valuable discussions during the preparation of
this paper. Also, many thanks to Alaeddin Malek for his supports
as the head of applied mathematics department of Tarbiat Modares
University.
\bibliographystyle{elsart-num-sort}
\bibliography{Majidian}
\end{document}